\documentclass{amsart}
\usepackage{amscd,amssymb,graphicx,epsfig}
\usepackage{array}

\title[When Polarizations Generate]{When Polarizations Generate}

\author{Gerald W. Schwarz}

 \address{Gerald W. Schwarz \newline
    \indent Department of Mathematics\newline
    \indent Brandeis University\newline
    \indent PO Box 549110\newline
    \indent Waltham, MA 02454-9110}
\email{schwarz@brandeis.edu}
\date{September 2006}

\newtheorem{theorem}[subsection]{Theorem}
\newtheorem{lemma}[subsection]{Lemma}
\newtheorem{proposition}[subsection]{Proposition}
\newtheorem{corollary}[subsection]{Corollary}
\newtheorem{conjecture}[subsection]{Conjecture}

\theoremstyle{definition}

\theoremstyle{remark} 
\newtheorem{remarks}[subsection]{Remarks}
\newtheorem{example}[subsection]{Example}

\newcommand{\op}{\operatorname}

\newcommand{\name}[1]{\textsc{#1\/}}
\newcommand{\NN}{{\mathbb N}} 

\newcommand{\ZZ}{{\mathbb Z}}

\newcommand{\CC}{{\mathbb C}}

\newcommand{\be}{\begin{enumerate}}
\newcommand{\ee}{\end{enumerate}}

\newcommand{\SL}{\op{SL}}

\newcommand{\GL}{\op{GL}}

\newcommand{\SO}{\op{SO}}

\newcommand{\Hom}{\op{Hom}}

\newcommand{\Orth}{\op{O}}

\newcommand{\inv}{^{-1}}

\newcommand{\quot}{/\!\!/}

\newcommand{\Sym}{\op{S}}

\renewcommand{\phi}{\varphi}
\newcommand{\pol}{\op{pol}}

 \numberwithin{equation}{subsection}

\thanks{Partially supported by  NSA Grant H98230-06-1-0023}
\begin{document}
\begin{abstract}
Let $G$ be a reductive complex  algebraic group and $V$ a finite-dimensional $G$-module. From elements of the invariant algebra $\CC[V]^G$ we obtain by polarization elements of  $\CC[kV]^G$, where $k\geq 1$ and $kV$ denotes the direct sum of $k$ copies of $V$.   For     $G$   simple our main result is the classification  of  the $G$-modules $V$ and integers $k\geq 2$ such that  polarizations generate $\CC[kV]^G$. 
\end{abstract}

\maketitle

\section{Introduction}\label{secintro}
Our base field is $\CC$, the field of complex numbers. Throughout this paper,  $G$ will denote a reductive algebraic group. All our $G$-modules are assumed to be  finite-dimensional and rational. Let $V$ be a $G$-module and let
 $f\in\CC[V]^G$ be homogeneous of degree $d$. For $v_1$, $v_2,\dots,v_k\in V$, consider the function $f(\sum_i s_i v_i)$ where the $s_i$ are indeterminates. Then
$$
f(\sum_i s_i v_i)=\bigoplus_{\alpha\in (\ZZ^+)^k,|\alpha|=d} s^\alpha f_\alpha(v_1,\dots,v_k) 
$$
where the $f_\alpha\in\CC[kV]^G$ are multihomogeneous of the indicated degrees $\alpha$.
Here  for $\alpha=(a_1,\dots,a_k)\in (\ZZ^+)^k$ we have $s^\alpha=s_1^{a_1} \cdots  s_k^{a_k}$ and $|\alpha|=a_1+\dots+a_k$. We call the $f_\alpha$ the \emph{polarizations} of $f$. Let $\pol_k(V)^G$ denote the subalgebra of $\CC[kV]^G$ generated by polarizations. 
We say that $V$ has the \emph{$k$-polarization property} if $\pol_k(V)^G=\CC[kV]^G$. 

 For $G=\Orth_n$ and $V$ the standard action on $\CC^n$  one has the $k$-polarization property for all $k$. We have the same result when   $G=S_n$ and $V$ is the standard action on $\CC^n$ by \name{Weyl} \cite{Weyl} and also for the standard actions of the Weyl groups of type $B$ and $C$ by   \name{Hunziker} \cite{Hunziker}. From  \name{Wallach} \cite{Wallach} we learn that the $2$-polarization property is false for the Weyl group of type $D_4$. 

One can also ask if $\CC[kV]^G$ is finite over $\pol_k[V]^G$. In case $V$ is a  module for $\SL_2$ which does not contain a copy of the irreducible two-dimensional module, then $\CC[kV]^G$ is finite over $\pol_k[V]^G$ for all $k$. See \name{Kraft-Wallach} \cite{KraftW} ($V$ irreducible) and \name{Losik-Michor-Popov} \cite{Popov} for the general case.  The same references prove the $k$-polarization property for all $k$ when $G$ is a torus. However, for the adjoint representation  of a simple group of rank at least $2$, finiteness fails for $k\geq 2$ \cite{Popov}.   See   \name{Draisma-Kemper-Wehlau} \cite{Wehlau} for questions about polarization and separation of orbits.

\begin{example}\label{sl2 example}
Let $V:=  \CC^2$  be the two-dimensional irreducible $G:=\SL_2$-module.  Then $\CC[V\oplus V]^G$ is generated by the determinant function. The polarizations of this generator give rise to three of the six determinant generators of $\CC[4 V]^G$. Hence $V\oplus V$ does not have the $2$-polarization property.
Of course, neither does $V$, since $\CC[V]^G=\CC$ while $\CC[V\oplus V]^G\neq\CC$.
\end{example}

Our main aim is to classify the $G$-modules which have the $k$-polarization property, $k\geq 2$, when $G$ is a simple linear algebraic group. Along the way we establish general criteria for a representation to have the $k$-polarization property. Note the trivial fact that the $(k+1)$-polarization property fails if the $k$-polarization property fails. The crucial case to consider is usually that of $k=2$.

We  thank V. Popov for useful comments.
\section{Slices}\label{secslices}
We establish some tools for obtaining our classification. Let $V$ be a $G$-module and
let $f\in\CC[V]^G$ be homogeneous of degree $d$. If  $k=2$ or $3$, then we will denote the polarizations of $f$ as $\{f_{i,j}\}_{i+j=d}$ and $\{f_{i,j,k}\}_{i+j+k=d}$, respectively.
We leave the proof of the following to the reader.
\begin{lemma}\label{lemsubmodule} Let $V$   be a $G$-module  with the $k$-polarization property. Then any $G$-submodule of $V$ has the $k$-polarization property.
\end{lemma}

 Let $V$ be a $G$-module. Then $\CC[V]^G$ is finitely generated, and we denote by $\pi\colon V\to V\quot G$ the morphism of affine varieties dual to the inclusion $\CC[V]^G\subset\CC[V]$. 
 
Let $V$ be a $G$-module and $v\in V$ such that the orbit $G\cdot v$ is closed. Then the isotropy group  $G_v$ is reductive and there is a splitting $V=S\oplus T_v(G\cdot v)$ of $G_v$-modules. The representation $G_v\to\GL(S)$ is called the \emph{slice representation at $v$}. We can arrange that $v\in S$.  We have a canonical map $\phi\colon G\ast^{G_v}S\to G\cdot S$ which is equivariant. Here $G\ast^{G_v}S$ is the quotient of $G\times S$ by the $G_v$-action sending
$(g,s)\mapsto (gh\inv,h\cdot s)$ for $h\in G_v$, $g\in G$ and $s\in S$. The $G_v$-orbit of $(g,s)$ is denoted $[g,s]$. Then $\phi\colon G\ast^{G_v}S\to V$ sends $[g,s]$ to $g(v+s)$. Replacing $S$ by an appropriate $G_v$-stable neighborhood  of $0\in S$ one has Luna's slice theorem \cite{Luna}. But here we only need one consequence of this theorem. Namely, that  the induced mapping $\phi\quot G\colon (G\ast^{G_v}S)\quot G\simeq S\quot G_v\to V\quot G$ induces an isomorphism of the Zariski cotangent spaces at the points $0$ and $G\cdot v$ in the quotients.   
 
\begin{lemma}\label{lemslice} Let $v$, $G_v$, etc.\ be as above. Suppose that $V$ has the $k$-polarization property. Then so does the   $G_v$-module  $S$.
\end{lemma} 
\begin{proof} We treat the case $k=2$ and leave the general case to the reader.
A $G_v$-stable complement to $T_{(v,0)}(G\cdot (v,0))$ in $V\oplus V$   is $S\oplus V$.  
Let $I:=\{f\in\CC[V\oplus V]^G\colon f(v,0)=0\}$, 
let $J:=\{f\in \CC[S\oplus V]^{G_v}\colon f(0,0)=0\}$ and 
let $K:=\{f\in\CC[S\oplus S]^{G_v}\colon f(0,0)=0\}$.   Let $\psi\colon  G\ast^{G_v}(S\oplus V)\to (V\oplus V)$ be the canonical map. Then, as indicated above, Luna's slice theorem implies that $\psi^*$ induces an isomorphism of $I/I^2$ with $J/J^2$, and clearly $J/J^2\to K/K^2$ is surjective. Thus $I\to K/K^2$ is surjective.
 
 Let $f\in\CC[V]^G$ be homogeneous of degree $d>0$. Then $f-f(v)\in I\cap \CC[V]^G$ and its image in $K\cap \CC[S]^{G_v}$ is the mapping  sending $s\in S$ to $f(v+s)-f(v)$. We have $f(v+s)-f(v)= \sum_{i+j=d, j>0} f_{i,j}(v,s)$. Thus each homogeneous component $f_{i,j}(v,s)$, $j>0$, lies in $K\cap \CC[S]^{G_v}$. If we polarize  in the second argument we obtain the collection of functions $\{ f_{i,j,k}(v,s,s')\colon j+k>0\}\subset K$ where $s$, $s'\in S$. On the other hand, if we take the polarizations $f_{d-k,k}$ of $f$ where $k>0$, then they are in $I$ and their images in $K$ are sums of the elements of $\{f_{i,j,k}(v,s,s')\colon k>0\}$. Thus the images of the $f_{i,j}-\delta_{0j}f(v)$ in $K$ are polarizations of elements of $K\cap\CC[S]^{G_v}$.  Since $V$ has the $2$-polarization property,   the $f_{i,j}-\delta_{0j}f(v)$  generate $I$ (as one varies $f$), hence $K/K^2$ is generated by the polarizations of elements of $K\cap\CC[S]^{G_v}$. But  functions in $K$ which span $K/K^2$ generate $\CC[S\oplus S]^{G_v}$, so $S$ has the $2$-polarization property.
 \end{proof}

 \begin{remarks} 
 \begin{enumerate}
 \item 
 The $G_v$-fixed part of $S$ plays no role. If $S'$ is the sum of the nontrivial isotypic components of $S$, then the interesting fact is that $\CC[S'\oplus S']^{G_v}$ is generated by polarizations. 
 \item Suppose that $\dim\CC[V]^G=1$ and the representation is stable (i.e., there is a non-empty open set of closed orbits). Then for any closed non-zero orbit $G\cdot v$, the slice representation is trivial, so that the Lemma is of no help.
\end{enumerate}
\end{remarks}
\section{Representations without the polarization property}

One can say that ``most'' representations do not have the $2$-polarization property. This is born out by the following sequence of lemmas.

 \begin{lemma} \label{lemnopolar}
 Suppose that  $V=V_1\oplus V_2$ where the $V_i$ are $G^0$-stable and the elements of $G$ preserve the $V_i$ or interchange them. For example, the $V_i$ could be the isotypic components corresponding  to  nontrivial irreducible $G^0$-modules.  Suppose further that $\CC[V_1\oplus V_2]^G$ has a minimal bihomogeneous generator $f$ of degree $(a,b)$ where $ab\geq 2$. Then $V$ does not have the $2$-polarization property.
 \end{lemma}
 \begin{proof} Let $v_i$, $v_i'$ denote elements in $V_i$, $i=1$, $2$.
 Set $d=a+b\geq 3$. Consider the polarization $f_{d-2,2}$ of $f$. We can write   $f_{d-2,2}(v_1,v_2,v_1',v_2')$ as a sum of terms $f^{2,0}+f^{1,1}+f^{0,2}$ where $f^{i,j}(v_1,v_2,v_1',v_2')$ has homogeneity $(i,j)$ in $v_1'$ and $v_2'$. The $f^{i,j}$ are $G'$-invariant, where $G'$ is the subgroup of $G$ preserving the $V_i$. Clearly, up to a scalar,   $G$ leaves $f^{2,0}+f^{0,2}$ and $f^{1,1}$ invariant. Hence the functions are $G$-invariant (since their sum is $G$-invariant) and nonzero (since $d\geq 3$).   Let $I$ denote the ideal of elements of $\CC[2V]^G$ vanishing at $0$. Suppose that   $\alpha f^{1,1}+\beta(f^{2,0}+f^{0,2})\in I^2$ for some $\alpha$ and $\beta$. We may assume that $f^{2,0}\neq 0$, i.e., that $a\geq 2$. Then evaluating at points $(v_1,v_2,v_1,0)$ we see that $\beta f\in I^2$, hence $\beta=0$. Now one evaluates at points $(v_1,v_2,v_1,v_2)$ to see that $\alpha=0$. We have shown that $f^{1,1}$ and $f^{2,0}+f^{0,2}$ are linearly independent modulo $I^2$.
  
  If $f^{1,1}$ is in the subalgebra generated by polarizations of   elements of $\CC[V]^G$, then   $f^{1,1}$ is a sum of terms $q  r_{\bullet-2,2}$ and $s  t_{\bullet-1,1}u_{\bullet-1,1}$ 
for appropriate homogeneous   $q,\dots,u\in\CC[V]^G$. Since $f^{1,1}$ is a minimal generator, our sum has to contain  terms of the form $r_{d-2,2}$.  Thus we may assume that $f^{1,1}\in r_{d-2,2}+I^2$ for some $r$. Restituting   we see that $(ab)f\in \binom d2 r+I^2\cap \CC[V]^G$.  Hence we have that $r+I^2=cf+I^2$ for some $c\neq 0$. It follows that $f^{1,1}\in c(f^{1,1}+f^{2,0}+f^{0,2})+I^2$ which implies that $f^{1,1}$ and $f^{2,0}+f^{0,2}$ are linearly dependent modulo $I^2$. This is a contradiction, hence $V$ does not have the $2$-polarization property.
 \end{proof}

  \begin{lemma} \label{lemsymplectic} Suppose that $G$ acts on $G^0$ by inner automorphisms and that  $V$ is a $G$-module which contains an irreducible symplectic $G^0$-submodule $U$. Further suppose that $\CC[U]^{G^0}$  has generators of even degree.  Then $V$ does not have the $2$-polarization property.
  \end{lemma}
  
  \begin{proof}    We may suppose that, as $G^0$-module, $V$ is the isotypic component of type $U$. A central torus of $G^0$ must act trivially on $U$, so we can reduce to the case that  $G^0$ is semisimple. Set $H:=Z_G(G^0)$. Then $H$ is finite and $G=HG^0$ where $H\cap G^0=Z(G^0)$.     Now $W:=\Hom(U,V)^{G^0}$ is an $H$-module (via the action of $H$ on $V$) and we have a canonical $G$-equivariant isomorphism of $V$ with $W\otimes U$ where the latter is naturally an $(H\times G^0)/Z(G^0)\simeq G$-module.    Let $k$ be minimal such that $\Sym^{2k}(W^*)^H\neq 0$.  We have a copy of  $\Sym^2(W^*)\subset \Sym^2(W^*)\otimes \wedge^2(U^*)^{G^0}\subset (W^*\otimes U^*\otimes W^*\otimes U^*)^{G^0}\subset \Sym^2(2V^*)^{G^0}$. Then there is a nonzero $f$ in the copy of $S^{2k}(W^*)^H$ in $\Sym^k(\Sym^2(W^*))$. 
   Suppose that $f$ is a polynomial in polarizations of generators of $\CC[V]^G$. We may assume that  each generator  lies in $\Sym^\lambda(W^*)^H\otimes\Sym^\lambda(U^*)^{G^0}$ for some partition $\lambda$ (so $\Sym^\lambda$ denotes the corresponding Schur component).  Then $f$ must be a polynomial in   polarizations of elements of the $\Sym^\lambda(W^*)^H\otimes\Sym^\lambda(U^*)^{G^0}$  where   $\Sym^\lambda$ is a symmetric power. But then by minimality of $k$ and the assumption that generators of $\CC[U]^{G^0}$ have even degree,  we must have that $f$ itself is a polarization. But, by construction, $f$ restitutes to $0$, i.e., $f(v,v)=0$, $v\in V$, hence $f$ is not a   polarization. 
   \end{proof}
  
  Let $R_j$ denote the irreducible $\SL_2$-module of dimension $j+1$, $j\in\NN$. 
  
  \begin{corollary}\label{corsl2} 
  Suppose that $G^0=\SL_2$ and that $V$ is a $G$-module which contains a $G^0$-submodule $R_j$ where $j$ is odd. Then $V$ does not have the $2$-polarization property.
  \end{corollary}
  
  \begin{proof} For $j$ odd, $R_j$ is a symplectic representation of $\SL_2$. Moreover, $\pm I \in\SL_2$   act  as $\pm 1$ on $R_j$, so that all elements of $\CC[R_j]^{\SL_2}$ have  even degree. Finally, all automorphisms of $\SL_2$ are inner. Thus we can apply Lemma~\ref{lemsymplectic}
  \end{proof}
  
  Now assume that $G^0=\CC^*$. Let $\nu_j$ denote the irreducible $\CC^*$-module with weight $j$. We denote by $ m\nu_j$ the direct sum of $m$ copies of $\nu_j$. Assume that   $V$ is a $G$-module such that the multiplicity of each $\nu_j$ is the same as that of $\nu_{-j}$ for all $j$. We say that $V$ is \emph{balanced} and we let $q(V)$ denote half the number of nonzero weight spaces, counting multiplicity.  
  
\begin{proposition} \label{proptorus}
Suppose that $G^0=\CC^*$. Let $V$ be a balanced $G$-module with $q(V)\geq 2$. Then $V$ does not have the $2$-polarization property.
\end{proposition}
 
\begin{proof}  First assume that there is a nonzero weight $j$ of multiplicity $m\geq 2$. Then the $\CC^*$-submodule $m(\nu_j\oplus\nu_{-j})$ is $G$-invariant, so that we may assume that it is all of $V$. We may then also assume that $j=1$.   
Set $V_1:=m\nu_1$ and $V_2:=m\nu_{-1}$. If $f$ is a minimal homogeneous generator of degree at least 3, then we are done by Lemma \ref{lemnopolar}. Thus
 we may  suppose that all the minimal homogeneous generators of $\CC[V_1\oplus V_2]^G$ have degree 2.  Let $G':=Z_G(G^0)$. Then $V_1$ and $V_2$ are $G'$-modules. Write $V_1=\oplus W_i$ where the $W_i$ are irreducible $G'$-submodules, and similarly write $V_2=\oplus U_ j$. Then the quadratic $G'$-invariants correspond to pairs $W_i$ and $U_j$ such that $U_j\simeq W_i^*$. But $\CC[V]^{G^0}$ has to be finite over $\CC[V]^{G'}$ and this forces that $V_1\simeq mW$ and $V_2\simeq mW^*$ for some irreducible one-dimensional representation $W$ of $G'$.  It follows that the image of $G'$ in $\GL(V)$ is that of $G^0=\CC^*$, so we may assume that $G'=G^0$. If $G=G^0$, then one can easily see that there are more quadratic generators in $\CC[2V]^G$ than those coming from polarizations. If $G\neq G^0$, then $G$ is generated by $G^0$ and an element $\alpha$ such that $\alpha t \alpha\inv=t\inv$ for $t\in G^0$ and $\alpha^2\in G^0$. Now $\alpha^2$ is fixed under conjugation by $\alpha$ so that $\alpha^2=\pm 1$.  For an appropriate basis $v_1,\dots v_m$ of $V_1$ and $w_1,\dots,w_m$ of $V_2$ we have that $\alpha(v_i)=w_i$ and $\alpha(w_i)=\pm v_i$, $i=1,\dots,m$. If $\alpha^2=-1$, then  one can see that  quadratic invariants do not generate $\CC[V]^G$. If $\alpha^2=1$, then we just have $m$-copies of the standard representation of $\Orth_2$. Since $m\geq 2$, one easily sees that there are more generators of $\CC[V\oplus V]^G$ than polarizations.

Now suppose that $V$ contains two different pairs of weights. Then we can assume that $V=m_1(\nu_p\oplus\nu_{-p})\oplus\ m_2(\nu_q\oplus\nu_{-q})$ as $\CC^*$-module where   $p$ and $q$ are relatively prime and $m_1$, $m_2\geq 1$. There is then clearly a bihomogeneous minimal $G$-invariant of degree $(a,b)$ where $ab\geq 2$, so that we can again apply Lemma \ref{lemnopolar}
\end{proof}

 If $V$ is a $G$-module  where $G^0$ is simple of rank 1 then we define $q(V)$ as before, relative to the action of a maximal torus.

\begin{corollary} \label{coreven}
Suppose that $G^0$ is simple of rank $1$ and that $V$ is a $G$-module with $q(V)\geq 3$.  Then $V$ does not have the $2$-polarization property.
\end{corollary}

\begin{proof} By Corollary \ref{corsl2} we may assume that, as $G^0$-module, 
$V$ is the direct sum of copies of $R_j$, $j$ even. Let $v\in V$ be a nonzero zero weight vector. Then the $G$-orbit though $w$ is closed with isotropy group a finite extension of $\CC^*$ and slice representation $V'$ where $q(V')\geq 2$. By Proposition \ref{proptorus}, $(V',G_v)$ does not have the $2$-polarization property, hence neither does $V$.
\end{proof}
 
\section{The Main Theorem}

Recall that a $G$-module is called \emph{coregular\/} if $\CC[V]^G$ is a regular $\CC$-algebra.

\begin{proposition} Let $V$ be an irreducible representation of the simple algebraic group $G$.
If $V$ is not coregular, then $V$ does not have the $2$-polarization property.
\end{proposition}
 
 \begin{proof} The   representations $R_j$ of $\SL_2$ which are not coregular have $q(R_j)\geq 3$, hence they do not have the $2$-polarization property by Corollary \ref{coreven}.
By \cite[Remark 5.2]{coreg} we know that if $V$ is not coregular and the rank of $G$ is at least 2, then one of the following occurs 
 \begin{enumerate}
  \item There is a closed orbit $G\cdot v$ such that $G_v$ has rank $1$. Let $G_v\to\GL(S)$ be the slice representation. Then $S$ is balanced. If $G_v^0$ is simple, then $q(S)\geq 4$ and if $G_v^0\simeq \CC^*$, then $q(S)\geq 2$.
 \item $V=\Sym^3(\CC^4)$ and $G=\SL_4$.
 \end{enumerate}
 By Proposition \ref{proptorus} and Corollary \ref{coreven}  any representation in (1) above does not have the $2$-polarization property. Thus the only remaining case is (2).  Here the laziest thing to do is to use the program LiE \cite{Lie} to compute some low degree invariants of one or two copies of $V$. The first generator of $\Sym^*(V^*)^G$ occurs in $\Sym^8(V^*)$ and the dimension of the fixed space is $1$. In $\Sym^2(V^*)\otimes\Sym^6(V^*)$ there is a two-dimensional space of invariants, so that $V$ does not have the $2$-polarization property.
  \end{proof}
  
  We would not have had to use LiE if the following could be established.
  
  \begin{conjecture}
  Let $H\subset\GL(V)$ where  $H$  is finite and not generated by reflections. Then $V$ does not have the $2$-polarization property.
  \end{conjecture}
  In the following we use the notation of \cite{coreg} for the simple groups and their representations. We list such representations as pairs $(V,G)$. 
  \begin{theorem}
Let $V$ be an irreducible nontrivial representation of the simple algebraic group $G$. If $V$ has the $2$-polarization property, then, up to (possibly outer) isomorphism, the pair $(V,G)$ is on the following list.
\begin{enumerate}
\item $(\phi_1,A_n)$, $n\geq 2$.
\item $(\phi_1^2,A_n)$, $n\geq 1$.
\item $(\phi_2,A_n)$, $n\geq 4$.
\item $(\phi_1,B_n)$, $n\geq 2$ and $(\phi_1,D_n)$, $n\geq 3$.
\item $(\phi_1,G_2)$.
\item $(\phi_3,B_3)$.
\item $(\phi_1,E_6)$.  
\end{enumerate}
\end{theorem}
\begin{proof} One can verify from \cite{coreg} that all the listed representations have the $2$-polarization property. We must rule out all other cases.  The list of coregular representations is due to \name{Kac-Popov-Vinberg} \cite{Vinberg}, see also \cite{coreg}. 

There are several easy ways to see that an irreducible coregular representation   fails to have the $2$-polarization property. \ One of the following can occur:
\begin{enumerate}
\item [(i)]$V$ is a symplectic representation of $G$.  
\item [(ii)] $\CC[V]^G$ is minimally generated by  homogeneous elements of degrees $m_1,\dots,m_k$ where $\sum_i (m_i+1)<\dim\CC[2V]^G$ or $\sum_i (m_i+1)=\dim\CC[2V]^G$ and $(2V,G)$ is not coregular.
\item [(iii)] $(V,G)$ is the adjoint representation of $G$ where $G$ has rank at least two. Here there is a slice representation whose effective part is the adjoint representation of $A_2$, $B_2$ or $G_2$. Then we can apply (ii).
 \end{enumerate}
Examples of (i) are the representations $(\phi_3,A_5)$, $(\phi_3,C_3)$ and $(\phi_5,B_5)$.  The representations   $(\phi_3,A_n$, $n=6$, $7)$  are examples of (ii) as is  
$(\phi_1,F_4)$.  Now we mention the remaining irreducible coregular representations  that can't be decided by the criteria above.
\begin{enumerate}
\item [(a)] $(\phi_8,D_8)$. Here $\CC[V]^G$ has generators in degree 2, 8, \dots while $\CC[2V]^G$ has three bihomogeneous invariants of degree $(2,2)$.
\item [(b)] $(\phi_3,A_8)$. Here $\CC[V]^G$ has generators in degrees 12, 18, \dots while $\CC[2V]^G$ has a bihomogeneous generator in degree (3,3).
\item [(c)] $(\phi_1^2, B_n)$, $n\geq 2$ or $(\phi_1^2,D_n)$, $n\geq 3$. Up to a trivial factor, these are just the representations of the groups $\SO_n$ on $\Sym^2(\CC^n)$, $n\geq 5$. Let $e_1,\dots e_n$ be the standard basis of $\CC^n$. The slice representation at the point $e_1^2+\dots+e_{n-3}^2$ is, up to trivial factors, the sum  $(\Sym^2(\CC^{n-3}),\SO_{n-3})\oplus(\Sym^2(\CC^3),\SO_3)$ and the latter representation does not have the $2$-polarization property by (ii).
\item [(d)] $(\phi_2,C_n)$, $n\geq 3$. In case that $n=3$, the representation fails to have the $2$-polarization property by (ii). For $n\geq 4$, the representations have a slice representation which contains a factor $(\phi_2,C_3)$, similarly to case (c).
\end{enumerate}
The remaining representations $(\phi_4,A_7)$ and $(\phi_4,C_4)$ are handled as in (a) and (b).
\end{proof}

One can now use the tables of \cite{coreg} to see which irreducible representations have the $k$-polarization property for $k\geq 3$.

\begin{corollary}
Let $G$ be simple and $V$ irreducible with the $k$-polarization property, $k\geq 3$. Then $(V,G)$ and $k$ are  on the following list.
\begin{enumerate}
\item $(\phi_1,A_n)$, $n\geq 2$ and $3\leq k<n+1$.
\item $(\phi_1,B_n)$, $n\geq 2$ and $3\leq k<2n+1$.
\item $(\phi_1,D_n)$, $n\geq 3$ and $3\leq k<2n$.
\item $(\phi_3,B_3)$ and $k=3$.
\end{enumerate}
\end{corollary}

It is also easy to determine the reducible representations with the $k$-polarization property for $k\geq 2$. One uses Lemma~\ref{lemnopolar}, the criterion (ii) above and the tables of \cite{coreg}.

\begin{corollary} Let $V$ be the direct sum of at least two irreducible nontrivial $G$-modules, where $G$ is simple. If $V$ has the $k$-polarization property, $k\geq 2$, then, up to isomorphism, $G=\SL_n$, $n\geq 5$ and  $V=j\CC^n$ where $2\leq j<n/k$.
\end{corollary}
\vskip1cm
 
\end{document}